# On Cantor's important proofs

W. Mueckenheim

University of Applied Sciences, Baumgartnerstrasse 16, D-86161 Augsburg, Germany

mueckenh@rz.fh-augsburg.de

______________________________________________

**Abstract.** It is shown that the pillars of transfinite set theory, namely the uncountability proofs, do not hold. (1) Cantor's first proof of the uncountability of the set of all real numbers does not apply to the set of irrational numbers alone, and, therefore, as it stands, supplies no distinction between the uncountable set of irrational numbers and the countable set of rational numbers. (2) As Cantor's second uncountability proof, his famous second diagonalization method, is an impossibility proof, a simple counter-example suffices to prove its failure. (3) The contradiction of any bijection between a set and its power set is a consequence of the impredicative definition involved. (4) In an appendix it is shown, by a less important proof of Cantor, how transfinite set theory can veil simple structures.

## 1. Introduction

Two finite sets have the same cardinality if there exists a one-to-one correspondence or bijection between them. Cantor extrapolated this theorem to include infinite sets as well: If between a set $M$ and the set $\mathbb{N}$ of all natural numbers a bijection $M \leftrightarrow \mathbb{N}$ can be established, then $M$ is denumerable or countable[1] and it has the same cardinality as $\mathbb{N}$, namely $\aleph_0$. There may be many non-bijective mappings, but at least *one* bijective mapping must exist. In 1874 Cantor [1] published the proof that the set $\mathbb{A}$ of all algebraic numbers (including the set $\mathbb{Q}$ of all rational numbers) is denumerable. His major achievement consists in having shown that the set $\mathbb{R}$ of all real numbers is uncountable, i.e. that *any* bijection $\mathbb{N} \leftrightarrow \mathbb{R}$ is impossible [1, 2]. He concluded that the cardinality of $\mathbb{R}$, denoted by $c$, surpasses $\aleph_0$, and he took this observation as the basis of transfinite set theory. Later on he simplified this proof by his celebrated second diagonalization method [3]. The third famous proof concerning the foundations of the concept of uncountability, due to Hessenberg [4], utilises the power set of the natural numbers.

It will be shown that the concept of uncountability fails because none of these proofs is conclusive.

______________________________

[1] The meaning of *countable set* covers *finite set* as well as *denumerable infinite set*.

## 2. Cantor's first proof of the uncountability of the real numbers

After long, hard work including several failures [5, p. 118 and p. 151] Cantor found his first proof showing that the set $\mathbb{R}$ of all real numbers cannot exist in form of a sequence. Here Cantor's original theorem and proof [1, 2] are sketched briefly, using his own symbols.

**Theorem 1.** *Consider an infinite sequence of different real numbers*

$$\omega_1, \omega_2, ..., \omega_\nu, ... \qquad (1)$$

*which is given by any rule, then we can find in any open interval ($\alpha$, $\beta$) a number $\eta$ (and, hence, infinitely many of such numbers) which is not a member of sequence (1).*

**Proof.** Take the first two members of sequence (1) which fit into the given interval ($\alpha$, $\beta$). They form the interval ($\alpha'$, $\beta'$). The first two members of sequence (1) which fit into this interval ($\alpha'$, $\beta'$) form the interval ($\alpha''$, $\beta''$) and so on. The result is a sequence of nested intervals. Now there are only two possible cases.

*Either* the number of intervals is finite. Inside the last one ($\alpha^{(\nu)}$, $\beta^{(\nu)}$) there cannot be more than one member of the sequence. Any other number of this interval ($\alpha^{(\nu)}$, $\beta^{(\nu)}$) can be taken as $\eta$.

*Or* the number of intervals is infinite. Then both, the strictly increasing sequence $\alpha, \alpha', \alpha'', ...$ and the strictly decreasing sequence $\beta, \beta', \beta'', ...$ converge to different limits $\alpha^\infty$ and $\beta^\infty$ or they converge to the same limit $\alpha^\infty = \beta^\infty$ (a case which always occurs in $\mathbb{R}$). $\alpha^\infty = \beta^\infty = \eta$ is not a member of sequence (1). If $\alpha^\infty < \beta^\infty$, then any number of [$\alpha^\infty$, $\beta^\infty$] satisfies theorem 1.

So far Cantor's proof. We note that the "*either*"-case cannot occur if sequence (1) contains at least all rational numbers (due to the denumerability of $\mathbb{Q}$ this is possible) because any interval ($\alpha^{(\nu)}$, $\beta^{(\nu)}$) $\subset \mathbb{Q}$ contains two rational numbers forming the next interval ($\alpha^{(\nu+1)}$, $\beta^{(\nu+1)}$) $\subset \mathbb{Q}$. Assuming this, we need only consider the "*or*"-case.

Cantor's proof is very general, but not without restrictions. In his first formulation quoted above, only sequences with *different* numbers are admitted. This restriction, however, is not really necessary and has later been abandoned [6]. But we note that, implicitly, there is a serious restriction with respect to the set of numbers it is to be applied to: the *complete* set of real numbers is required as the underlying manifold investigated. If only one of them is removed, the proof fails because just this one could be the common limit $\alpha^\infty = \beta^\infty$.

Cantor took his result as evidence in favour of the existence and uncountability of the set $\mathbb{T}$ of all transcendental numbers which were shortly before discovered by Liouville [7]. Nevertheless his proof fails, if applied to the set $\mathbb{T}$ alone. The reason is again, that "any infinite sequence" like $\alpha, \alpha', \alpha'', ...$ or $\beta, \beta', \beta'', ...$ need not converge to a transcendental limit. Already the absence of a single number, zero for instance, cannot be tolerated, because it is the limit of several sequences.

This situation, however, is the same if only the set $\mathbb{Q}$ of all rational numbers is considered. Therefore both sets, $\mathbb{Q}$ and $\mathbb{T}$, have the same status with respect to *this* uncountability proof. And we are not able, *based on this very proof*, to distinguish between them.

On the other hand, the proof can feign the uncountability of a countable set. If, for instance, the alternating harmonic sequence

$$\omega_\nu = (-1)^\nu / \nu \to 0$$

is taken as sequence (1), yielding the intervals $(\alpha, \beta) = (-1, 1/2)$, $(\alpha', \beta') = (-1/3, 1/4)$, ... we find that its limit 0 does not belong to the sequence, although the set of numbers involved, $\mathbb{N} \cup \{0\}$, is obviously denumerable, e.g., by the bijection $0 \leftrightarrow 1$ and $\omega_\nu \leftrightarrow (\nu + 1)$ for $\nu \in \mathbb{N}$.

The alternating harmonic sequence does not, of course, contain all real numbers, but this simple example demonstrates that Cantor's first proof is not conclusive. *Based upon this proof alone*, the uncountability of this and every other alternating convergent sequence must be claimed. Only from *some other information* we know their countability (as well as that of $\mathbb{Q}$), but how can we exclude that *some other information*, not yet available, in the future will show the countability of $\mathbb{T}$ or $\mathbb{R}$?

Anyhow, the countability properties of an infinite set will not be altered by adding or removing one single element. The fact that Cantor's first uncountability proof does not apply to the set $\mathbb{R} \setminus \{r\}$, with $r$ being any real number, shows the failure of this proof if not the failure of the whole concept of countability.

## 3. Cantor's second diagonalization method

The first uncountability proof was later on [3] replaced by a proof which has become famous as Cantor's second diagonalization method (SDM). Try to set up a bijection between all natural numbers $n \in \mathbb{N}$ and all real numbers $r \in [0,1)$. For instance, put all the real numbers at random in a list with enumerated rows (by list we will understand an injective sequence, i.e., an injective function with domain $\mathbb{N}$). The digits of the $n^{\text{th}}$ real number may be denoted by

$$r_n = 0.a_{n1}a_{n2}a_{n3}... \tag{2}$$

**Table 1.** Attempt to list the set of all real numbers of [0, 1) (with underlined digits to be replaced)

| $n$ | $r_n$ |
|---|---|
| 1 | $0.\underline{a}_{11}a_{12}a_{13}...$ |
| 2 | $0.a_{21}\underline{a}_{22}a_{23}...$ |
| 3 | $0.a_{31}a_{32}\underline{a}_{33}...$ |
| ... | ... |

If the diagonal digit $a_{nn}$ of each real number $r_n$ is replaced by

$$b_{nn} \neq a_{nn} \quad \text{with} \quad 1 \leq b_{nn} \leq 8 \tag{3}$$

(in order to avoid identities like 0.999... = 1.000...) we can construct the diagonal number

$$R = 0.b_{11}b_{22}b_{33}... \tag{4}$$

belonging to the real interval [0, 1) but differing from any $r_n$ of the list. This shows that the list of table 1 and any such list must be incomplete.

But the impossibility of any bijection between infinite sets cannot be shown by some examples; it must be based upon a general contradiction. This means, Cantor's SDM must show that, *in any conceivable case*, a number $R \neq r_n$ is constructed, which was not on the list; no restriction must be imposed concerning its structure and contents; no rule other than (3) must restrict the replacement of the diagonal digits to construct $R$ according to (4). Otherwise, the possibility of a bijection is not excluded. These conditions are observed in Cantor's original paper [3] and also in all text books reporting on this topic.

Now we will contradict Cantor's SDM by a special counter-example. We will construct a list which belongs to the set of allowed lists (in principle, it could have been generated by accident), and which is easy to treat by an allowed replacement process (in principle, it could happen by accident). But the result is not a new diagonal number but a diagonal number contained in the list.

**Table 2.** Special list with real numbers (with underlined digits to be replaced by 1)

| $n$ | $r_n$ |
|---|---|
| 1 | 0.0̲00... |
| 2 | 0.10̲00... |
| 3 | 0.110̲00... |
| ... | ... |

This list is a rational sequence with limit 1/9. If now the diagonal number is constructed, always replacing the diagonal digit 0 by 1, we obtain, after changing $n$ (> 6) digits, the diagonal number

$$R(n) = 0.111...111 = r_{n+1} \qquad (5)$$

with $n$ digits 1, equal to the real number of line number $n + 1$ and, hence, being *not* different from any number of the list, not even from any rational number. This situation does not change with growing $n$ approaching infinity. And more cannot be done.

Further, this second uncountability proof suffers from the same handicap as the first one (and as all the proofs to be discussed): A general contradiction is claimed but when removing only one element of the interval [0, 1), then the general contradiction cannot be obtained.

## 4. The set of functions

Another application of SDM is to show that the set $F$ of all real functions $f(x)$ in the interval [0, 1) has a cardinality larger than c, that of the continuum. The proof runs in short as follows [8, 9]. Try to set up a bijection between the interval [0, 1) and the set $F$ of all real functions $f(x)$

$$y \leftrightarrow f_y(x) \tag{6}$$

where $y \in [0, 1)$. Consider the function $f_y(x)$, related to $y$, calculate its value at position $x = y$, namely $f_y(y)$, and increase it by 1, resulting in $f_y(y) + 1$. Combine all the values obtained in that way to create the function $g(y) = f_y(y) + 1$. This is a function in the interval $[0, 1)$ which differs from any $f_y(x) \in F$ at least in one point, namely for $x = y$. So (6) is not a bijection; the set $F$ has a cardinality $f > c$.

This proof can be contradicted by the same method as applied in section 3. As the attempted bijection (6) is in no way specified, we can require that the function related to $(1 + y)/2 \in [0, 1)$

$$f_{(1 + y)/2}(x) = g(x) \tag{7}$$

is always the same as the diagonal function $g(y)$ constructed up to the point $y$. Then no new function is created at all. This impossibility proof is void too.

## 5. The power set of $\mathbb{N}$

A general proof by Hessenberg [4] shows that there is no bijection between $\mathbb{N}$ and its power set $\mathcal{P}(\mathbb{N})$. If $\mathbb{N}$ could be bijected with its power set $\mathcal{P}(\mathbb{N})$, then some $n \in \mathbb{N}$ unavoidably would be mapped on subsets $s(n)$ not containing them. The subset $M$ of all such numbers $n$

$$M = \{n \mid n \rightarrow s(n) \text{ and } n \notin s(n)\} \tag{8}$$

belongs to $\mathcal{P}(\mathbb{N})$ as an element. But the set $M$ together with the mapping $s$ does not exist. If $M$ does not contain the element $m$ which is mapped upon it by $s: m \rightarrow M$, then $m$ belongs to $M$, but exactly then $M$ must not contain $m$ and so on: $m \in M \Rightarrow m \notin M \Rightarrow m \in M \Rightarrow m \notin M \Rightarrow ...$

This paradox is commonly interpreted as contradicting the existence of any bijection $\mathbb{N} \leftrightarrow \mathcal{P}(\mathbb{N})$. But it is clear that the set $\{M, m, s\}$ belongs to the class of impredicatively defined sets like Russell's set of all sets which do not contain themselves. In order to determine $M$, the mapping $s$ must already be complete. But $s$ is defined by the condition that the pre-image $m$ of $M$ is not an element of $M$ if it is an element of $M$. Impossible sets like the set of all sets are well known for their power of generating paradoxes and have been banned from set theory long ago.

The triple {*M*, *m*, *s*} has not yet been recognized as an impossible set because each of its members on its own can exist. But this fact must not be misinterpreted as an uncountability proof of the powerset of $\mathbb{N}$. For an example of two equivalent sets which do not allow for a mapping involving the triple {*M*, *m*, *s*} see [10].

### 6. Concluding remarks

The three proofs discussed above (together with some variations) are the only proofs from which evidence for the existence of transfinite cardinal numbers can be obtained. As we have seen, in fact they prove nothing at all. This leads to the result that the assumption of countable and uncountable sets and the assumption of different infinities are invalid. Note that none of the concepts and results of transfinite set theory have found use in any other science.

There is no actual infinity, because it would lead to contradictions. Cantor's continuum hypothesis, i.e., the question whether c is the next higher transfinite number $\aleph_1 = 2^{\aleph_0}$, is not only undecidable but meaningless. There is only the well known (more or less) potential infinity. So we can exorcise all the alephs, beths and higher accessible and (hyper-) inaccessible cardinal numbers and be satisfied with the conventional symbol $\infty$ introduced by Wallis for the only existing potential infinity. And the equations

$$\infty = 1 + \infty = 2\infty = \infty\infty \tag{9}$$

can be supplemented by

$$\infty = 2^\infty = \infty^{\infty^{\infty^{\cdots}}}. \tag{10}$$

Nowhere and never does infinity become actual or complete or finished. It has no value. It merely marks a direction and, therefore, does tolerate nearly every equation like (9) and (10).

# Appendix.  On Cantor's proof of continuity-preserving manifolds

A less important but very instructive proof of Cantor [6] is analysed below, which shows in a striking manner how the use of transfinite set theory veils even most simple structures.

## A1. Introduction

The manifold $\mathbb{R}^n$ (with $n \geq 2$) remains continuous if the set of points with purely algebric co-ordinates is taken off. According to Cantor's interpretation this is *a peculiar property of countable sets*.[2] In fact, this property does not only hold for all countable sets but it is the same for many uncountable sets. It becomes immediately clear, by simplifying the proof, that the continuity of $\mathbb{R}^n$ is also preserved, if the uncountable set of points with purely transcendental co-ordinates is removed.

In the following the expressions "manifold" and "continuum" denote the $n$-dimensional Euclidean space. They are used synonymously with the set of points, each of which is determined by a set of $n$ co-ordinates, which the $n$-dimensional Euclidean space is isomorphic to.

By defining origin and axes of a co-ordinate system the points of a manifold are subdivided into three sets: the countable set $\mathbb{AA}$ of those points with purely algebric numbers as co-ordinates, the uncountable set $\mathbb{TT}$ of those points with purely transcendental numbers as co-ordinates, and the uncountable set $\mathbb{AT}$ of the remaining points with mixed co-ordinates, i. e., with at least one algebric number and at least one transcendental number serving as a co-ordinate. Of course these properties do not belong to a point itself because the type of co-ordinate system as well as its origin and its axes can be chosen in an arbitrary way. But once the system has been fixed, the bijective mapping of the points $N$ of the continuum on the $n$-tuples $(x_1, x_2, ..., x_n)$ is fixed too

$$N \leftrightarrow (x_1, x_2, ..., x_n) \in \mathbb{R}^n = \mathbb{AA} \cup \mathbb{TT} \cup \mathbb{AT}. \tag{A1}$$

---

[2] Was die abzählbaren Punktmengen betrifft, so bieten sie eine merkwürdige Erscheinung dar, welche ich im folgenden zum Ausdruck bringen möchte. [5, p. 154]

A line or curve $l$ connecting two points of $\mathbb{R}^n$ may contain infinitely many points of $\mathbb{A}\!\mathbb{A}$. If the latter set is taken off, $l$ is no longer continuous in the remaining manifold ($\mathbb{R}^n \setminus \mathbb{A}\!\mathbb{A}$). But in 1882 Cantor [6] proved that $\mathbb{R}^n$ (with $n \geq 2$) remains continuous even if the set of points with purely algebraic co-ordinates is taken off. This means between two of the remaining points, with not purely algebraic co-ordinates, which Cantor called $N$ and $N'$, one can always find a continuous linear connection of the same character, which Cantor called $l'$. In short

$$N, N' \in (\mathbb{R}^n \setminus \mathbb{A}\!\mathbb{A}) \quad \Rightarrow \quad \exists\, l'(N, N'):\ N, N' \in l' \subset (\mathbb{R}^n \setminus \mathbb{A}\!\mathbb{A}) \qquad (A2)$$

## A2. Cantor's proof of the continuity of $(\mathbb{R}^n \setminus \mathbb{A}\!\mathbb{A})$

The set $\mathbb{A}\!\mathbb{A}$ is countable. Hence, any interval of the uncountable set $l$ contains points belonging to the uncountable set $(\mathbb{R}^n \setminus \mathbb{A}\!\mathbb{A})$. We consider a finite set of them $\{N_1, N_2, ..., N_k\}$. Between any pair of these points a part of a circle can be found which connects these points but contains no point of $\mathbb{A}\!\mathbb{A}$. This is shown for two points, $N$ and $N_1$, as follows: The centres of circles which on their circumference contain at least one point of $\mathbb{A}\!\mathbb{A}$ form a countable set. The centres of circles containing on their circumference $N$ and $N_1$ belong to a straight line (i.e., an uncountable set). This line contains at least one point which is centre of a circle containing on its circumference $N$ and $N_1$ but not any point of $\mathbb{A}\!\mathbb{A}$. As this can be shown for any pair of the finite set of points $\{N, N_1, N_2, ..., N_k, N'\}$ the proof is complete (see Fig. 1 (a) for the two-dimensional case).

The proof would work as well, if only the original pair of points, $N$ and $N'$, had been considered.

## A3. Simplified proof of the continuity of $(\mathbb{R}^n \setminus \mathbb{A}\!\mathbb{A})$

We apply Cartesian co-ordinates. At least one co-ordinate, say $x_\nu$, of the not purely algebraic point

$$N = (x_1, ..., x_{\nu-1}, x_\nu, x_{\nu+1}, ..., x_n) \qquad (A3)$$

is transcendental. Let this be constant while all the other $x_\mu$ (with $\mu = 1, ..., \nu-1, \nu+1, ..., n$) are continuously changed until they reach the values of the co-ordinates of $N'$

$$(x'_1, ..., x'_{\nu-1}, x_\nu, x'_{\nu+1}, ..., x'_n) \qquad (A4)$$

The $x_\mu$ (with $\mu = 1, ..., \nu-1, \nu+1, ..., n$) define a hyper plane $\mathbb{R}^{n-1} \subset \mathbb{R}^n$, within which we can choose an arbitrary way. If at least one of the final co-ordinates $x'_\mu$ is transcendental, we finish the proof by changing $x_\nu$ to $x'_\nu$ without leaving the set $\mathbb{T} \cup \mathbb{AT}$ (see Fig. 1 (b) for the two-dimensional case). If none of the final co-ordinates $x'_\mu$ is transcendental, we stop the process of continuously changing the $x_\mu$ for one of those co-ordinates, $x_\rho$, at the transcendental value $x''_\rho$ before the final algebraic value $x'_\rho$ of (A4) is reached (or we re-adjust $x''_\rho$ afterwards). Then, staying always in $\mathbb{T} \cup \mathbb{AT}$ we let $x_\nu$ approach $x'_\nu$, which in this case must be transcendental, and finally we complete the process by changing $x_\rho$ from its intermediate transcendental value $x''_\rho$ to its final value $x'_\rho$ (see Fig. 1 (d) for the two-dimensional case).

This method can also be used between any pair of points $\{N, N_1, N_2, ..., N_k, N'\}$ of $l$ belonging to the set $\mathbb{T} \cup \mathbb{AT}$. Though the complete length of the connection $l'$ remains unchanged, the deviation of any of its points from the straight line $l$ can be made as small as desired (see Fig. 1 (c) for the two-dimensional case).

## A4. Proof of the continuity of $(\mathbb{R}^n \setminus \mathbb{T})$

Cantor, who considered the preserved continuity of $(\mathbb{R}^n \setminus \mathbb{AA})$ a peculiar property of countable sets, obviously overlooked that taking off the uncountable set $\mathbb{T}$ leads to a continuous manifold too, similar to that remaining after taking off $\mathbb{AA}$. This fact becomes immediately clear from the proof given in section A3 but remains veiled in Cantor's original version given in section A2.

We see immediately that any pair of points with at least one algebraic co-ordinate can be connected by a continuous linear subset of the same character

$$N, N' \in (\mathbb{R}^n \setminus \mathbb{T}) \quad \Rightarrow \quad \exists\, l'(N, N'): N, N' \in l' \subset (\mathbb{R}^n \setminus \mathbb{T}). \qquad (A5)$$

The proof runs precisely as demonstrated in section A3, with the only difference that those co-ordinates which there and in the caption of Fig. 1 are prescribed as transcendental, now have to be algebraic. We can even go further and take off all points with purely non-rational co-ordinates or even all points with purely non-natural co-ordinates, so that there remains at least one $x_\nu \in \mathbb{N}$ of point $N$ and at least one $x'_\mu \in \mathbb{N}$ of point $N'$. It is obvious then that $N$ and $N'$ have a continuous connection as depicted in Fig. 1 (b) or in Fig. 1(d) along the "grid lines". In fact there are infinitely many of these connections.

As an example we consider the points $N = (n, \xi)$ and $N' = (n', \xi')$ with $n, n' \in \mathbb{N}$ and $\xi, \xi' \in \mathbb{R}$ in the two-dimensional Cartesian co-ordinate system $\mathbb{R}^2$. After taking off all points except those with at least one co-ordinate being a natural number, we have in the remaining manifold the connection by changing co-ordinates as described in (A6). First, moving along the grid line $x_1 = n$, change $\xi$, the possibly non-natural $x_2$-co-ordinate of $N$, to an intermediate co-ordinate $x''_2 = m$, choosing any natural number $m$. Then change $x_1 = n$ to $x'_1 = n'$, moving along the grid line $x''_2 = m$. Finally change $x''_2 = m$ to $x'_2 = \xi'$, moving along the grid line $x'_1 = n'$, briefly

$$N = (x_1, x_2) = (n, \xi) \to (n, m) \to (n', m) \to (n', \xi') = (x'_1, x'_2) = N'. \qquad (A6)$$

The connection does not contain any point with purely non-natural co-ordinates.

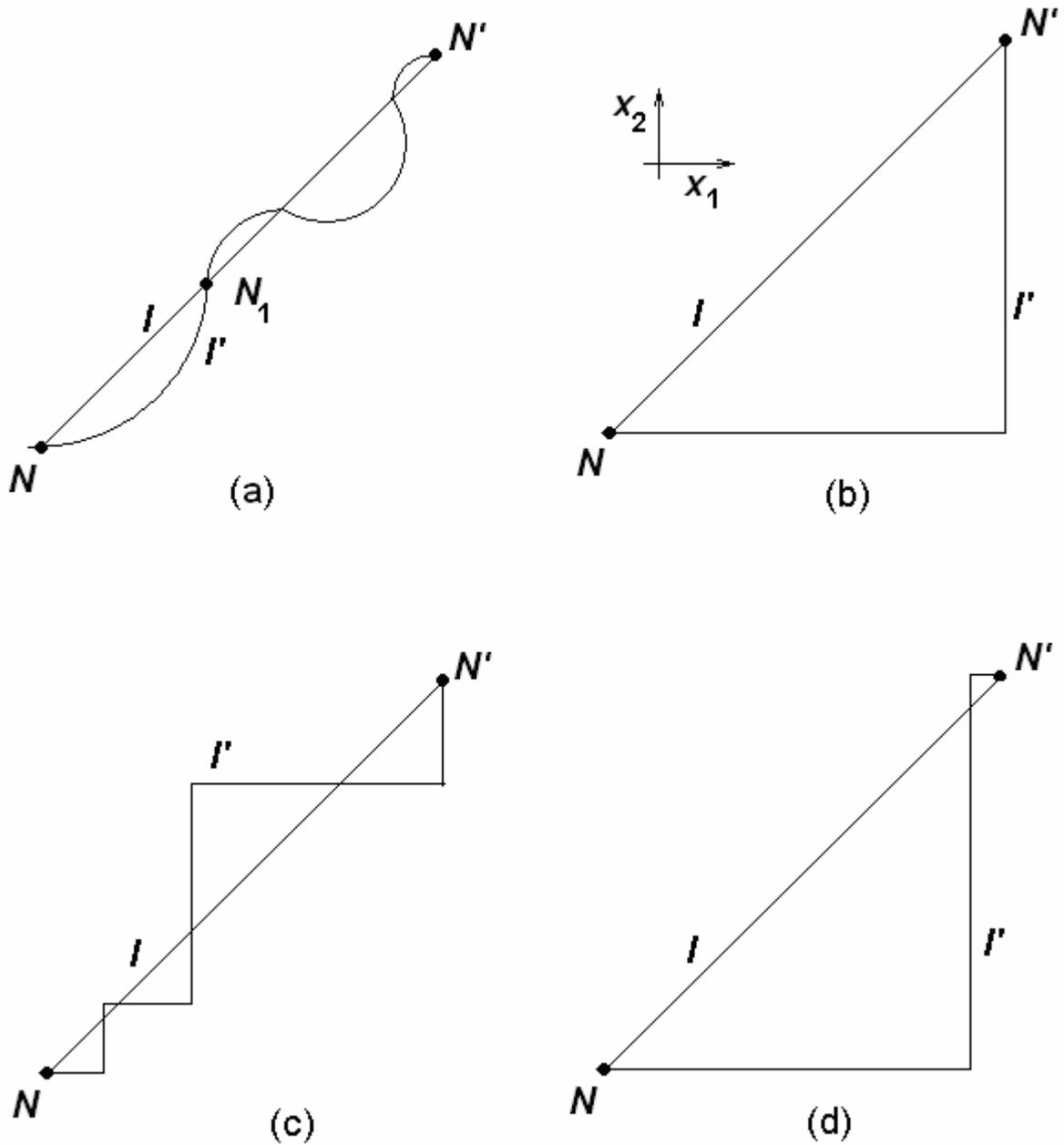

**Fig. 1. (a)** Connection $l'$ circumventing points of $\mathbb{A}\mathbb{A}$ between $N = (x_1, x_2)$ and $N' = (x'_1, x'_2)$ proposed by Cantor [6], **(b)** $l'$ of the present proof in case $x_2$ and $x'_1$ are transcendental, **(c)** same as (b) in case a smaller deviation of $l'$ from the straight line $l$ is requested, **(d)** $l'$ of the present proof in case only $x_2$ and $x'_2$ being transcendental.